\documentclass{amsart}%
\usepackage{amssymb}
\usepackage{amsmath}
\usepackage{latexsym}
\usepackage{hyperref}
\usepackage{amsfonts}
\usepackage{graphicx}%
\newcommand{\be}{\begin{equation}}
\newcommand{\ee}{\end{equation}}
\newcommand{\bea}{\begin{eqnarray}}
\newcommand{\eea}{\end{eqnarray}}

\begin{document}
\title{Mild Ricci curvature restrictions for steady gradient Ricci solitons}
\author{Bennett Chow}
\author{Peng Lu$^{1}$}
\maketitle

All objects are assumed to be $C^{\infty}$.\footnotetext[1]{Addresses. Bennett
Chow: Math. Dept., UC San Diego; Peng Lu: Math. Dept., U of Oregon.}\ Let
$\left(  \mathcal{M}^{n}\!,g\right)  $ be a complete Riemannian mani\-fold and
$\phi\!:\!\mathcal{M}\rightarrow(0,\infty)$.\ For $\gamma\!:\!\left[
0,\bar{s}\right]  \rightarrow\mathcal{M}$, $\bar{s}>0$, define (see Li and Yau
\cite{LY})\vspace{-0.02in}%
\[
\mathcal{J}\left(  \gamma\right)  =\int_{0}^{\bar{s}}\left(  \left\vert
\gamma^{\prime}\left(  s\right)  \right\vert ^{2}+2\phi\left(  \gamma\left(
s\right)  \right)  \right)  ds.\vspace{-0.02in}%
\]
For any $x,y\in\mathcal{M}$, among all paths joining $x$ to $y$, there exists
a minimizer of $\mathcal{J}$.

Given $\gamma_{u}:\left[  0,\bar{s}\right]  \rightarrow\mathcal{M}$,
$u\in\left(  -\varepsilon,\varepsilon\right)  $, define $\Gamma:\left[
0,\bar{s}\right]  \times\left(  -\varepsilon,\varepsilon\right)
\rightarrow\mathcal{M}$ by $\Gamma\left(  s,u\right)  =\gamma_{u}\left(
s\right)  $ and define $S=\Gamma_{\ast}\left(  \partial/\partial s\right)  $
and $U=\Gamma_{\ast}\left(  \partial/\partial u\right)  $. The first variation
formula is\vspace{-0.02in}%
\[
\tfrac{1}{2}\!\left.  \tfrac{d}{du}\right\vert _{u=0}\mathcal{J}\left(
\gamma_{u}\right)  =\int_{0}^{\bar{s}}\left(  \left\langle \nabla
_{S}U,S\right\rangle +U\left(  \phi\right)  \right)  ds=\int_{0}^{\bar{s}%
}\left\langle U,-\nabla_{S}S+\nabla\phi\right\rangle ds+\left.  \left\langle
U,S\right\rangle \right\vert _{0}^{\bar{s}}.\vspace{-0.02in}%
\]
Hence the critical points $\gamma$ of $\mathcal{J}$ on paths with fixed
endpoints, called $\phi$-geodesics, are the paths which satisfy $\nabla
_{S}S=\nabla\phi$. This implies $\left\vert S\right\vert ^{2}-2\phi
=C=\operatorname{const}$ on $\gamma$.

Let $\operatorname{Rm}$ denote the Riemann curvature tensor. The second
variation formula is\vspace{-0.02in}%
\begin{align*}
\tfrac{1}{2}\!\left.  \tfrac{d^{2}}{du^{2}}\right\vert _{u=0}\mathcal{J}%
\left(  \gamma_{u}\right)   &  =\int_{0}^{\bar{s}}\left(  \left\vert
\nabla_{S}U\right\vert ^{2}-\left\langle \operatorname*{Rm}\left(  S,U\right)
U,S\right\rangle +\nabla\nabla\phi\left(  U,U\right)  \right)  ds\\
&  \quad\;-\int_{0}^{\bar{s}}\left\langle \nabla_{U}U,\nabla_{S}S-\nabla
\phi\right\rangle ds+\left.  \left\langle \nabla_{U}U,S\right\rangle
\right\vert _{0}^{\bar{s}},\vspace{-0.02in}%
\end{align*}
where we used $\nabla\nabla\phi\left(  U,U\right)  =U\left(  U\left(
\phi\right)  \right)  -\left(  \nabla_{U}U\right)  \left(  \phi\right)  $. So
when $\gamma_{0}$ is a $\phi$-geodesic,\vspace{-0.02in}%
\[
\tfrac{1}{2}\!\left.  \tfrac{d^{2}}{du^{2}}\right\vert _{u=0}\mathcal{J}%
\left(  \gamma_{u}\right)  -\left.  \left\langle \nabla_{U}U,S\right\rangle
\right\vert _{0}^{\bar{s}}=\!\int_{0}^{\bar{s}}\!\left(  \left\vert \nabla
_{S}U\right\vert ^{2}-\left\langle \operatorname*{Rm}\left(  S,U\right)
U,S\right\rangle +\nabla\nabla\phi\left(  U,U\right)  \!\right)
\!ds.\vspace{-0.02in}%
\]
Now assume $\zeta:\left[  0,\bar{s}\right]  \rightarrow\mathbb{R}$ vanishes at
the endpoints, let $\left\{  e_{i}\right\}  _{i=1}^{n}$ be a parallel
orthonormal frame along $\gamma_{0}$, take $U=\zeta e_{i}$, and sum. If
$\gamma=\gamma_{0}$ is $\phi$-minimal, then\vspace{-0.02in}%
\[
0\leq\tfrac{1}{2}%
{\textstyle\sum\limits_{i=1}^{n}}
\left.  \tfrac{d^{2}}{du_{i}^{2}}\right\vert _{u_{i}=0}\mathcal{J}\left(
\gamma_{u_{i}}\right)  =\int_{0}^{\bar{s}}\left(  n\left(  \zeta^{\prime
}\right)  ^{2}+\zeta^{2}\Delta\phi-\zeta^{2}\operatorname{Rc}\left(
S,S\right)  \right)  ds.\vspace{-0.02in}%
\]
Let $f:\mathcal{M}\rightarrow\mathbb{R}$. By $\nabla\nabla f\left(
S,S\right)  =S\left(  S\left(  f\right)  \right)  -\left\langle \nabla
_{S}S,\nabla f\right\rangle $, and $\nabla_{S}S=\nabla\phi$, we have\vspace
{-0.02in}%
\[
\int_{0}^{\bar{s}}\zeta^{2}\operatorname{Rc}\left(  S,S\right)  ds=\int
_{0}^{\bar{s}}\zeta^{2}\operatorname{Rc}_{f}\left(  S,S\right)  ds+2\int
_{0}^{\bar{s}}\zeta\zeta^{\prime}\left\langle \nabla f,S\right\rangle
ds+\int_{0}^{\bar{s}}\zeta^{2}\left\langle \nabla\phi,\nabla f\right\rangle
ds,\vspace{-0.02in}%
\]
where $\operatorname{Rc}_{f}=\operatorname{Rc}+\nabla\nabla f$ and we
integrated by parts. Let $\Delta_{f}=\Delta-\nabla f\cdot\nabla$.
Then\vspace{-0.02in}%
\begin{equation}
-\int_{0}^{\bar{s}}\zeta^{2}\Delta_{f}\phi ds+\int_{0}^{\bar{s}}\zeta
^{2}\operatorname{Rc}_{f}\left(  S,S\right)  ds\leq\int_{0}^{\bar{s}}\left(
n\left(  \zeta^{\prime}\right)  ^{2}-2\zeta\zeta^{\prime}\left\langle \nabla
f,S\right\rangle \right)  ds.\vspace{-0.02in}\label{Exeter NH}%
\end{equation}

Let $\left(  \mathcal{M},g,f\right)  $ be a complete steady gradient Ricci
soliton with $\operatorname{Rc}_{f}=0$, $R+\left\vert \nabla f\right\vert
^{2}=1$ and $R>0$. Let $c>0$ and $\phi=cR$.\footnote[2]{So $\mathcal{J}$ is
like Perelman's $\mathcal{L}$-length \cite{Perelman1}. In contrast, the
integral curves to $\nabla f$ are $-\frac{R}{2}$-geodesics.} Since
$-\Delta_{f}R=2\left\vert \operatorname{Rc}\right\vert ^{2}$, $\left\vert
\nabla f\right\vert \leq1$, and $\left\vert \gamma^{\prime}\right\vert
^{2}=\left\vert S\right\vert \leq\sqrt{C+2c}$, on a minimal $cR$-geodesic we
have\vspace{-0.02in}%
\begin{equation}
\int_{0}^{\bar{s}}\zeta^{2}\left\vert \operatorname{Rc}\right\vert ^{2}%
ds\leq\frac{n}{2c}\int_{0}^{\bar{s}}\left(  \zeta^{\prime}\right)
^{2}ds+\frac{\sqrt{C+2c}}{c}\int_{0}^{\bar{s}}\left\vert \zeta\zeta^{\prime
}\right\vert ds.\vspace{-0.02in}\label{369 894}%
\end{equation}
Since $\left\vert \nabla R\right\vert =2\left\vert \operatorname{Rc}\left(
\nabla f\right)  \right\vert \leq2\left\vert \operatorname{Rc}\right\vert $,
we have the same estimate for $\frac{1}{4}\int_{0}^{\bar{s}}\zeta
^{2}\left\vert \nabla R\right\vert ^{2}ds$.

For $x,y\in\mathcal{M}$ with $d(x,y)\geq4$, let $\gamma\!:\!\left[  0,\bar
{s}\right]  \rightarrow\!\mathcal{M}$, $\bar{s}=d\left(  x,y\right)  $, be a
minimal $cR$- geodesic from $x$ to $y$, and let $\bar{\gamma}\!:\![0,\bar
{s}]\!\rightarrow\!\mathcal{M}$ be a minimal geodesic from $x$ to $y$.
Then\vspace{-0.02in}%
\begin{equation}
C\bar{s}  \leq\!\int_{0}^{\bar{s}}\!\!\left(  \left\vert \gamma^{\prime
}\left(  s\right)  \right\vert ^{2}+2cR(\gamma(s))\right)  \!ds\leq\!\int
_{0}^{\bar{s}}\!\!\left(  \left\vert \bar{\gamma}^{\prime}\left(  s\right)
\right\vert ^{2}+2cR(\bar{\gamma}(s))\right)  \!ds\leq(1+2c)\bar
{s},\label{C bounds}
\end{equation}
Hence $C\leq 1+2c$. If $\zeta\left(  s\right)  =s$ for $s\in\left[
0,1\right]  $, $\zeta\left(  s\right)  =1$ for $s\in\left[  1,\bar
{s}-1\right]  $, and $\zeta\left(  s\right)  =\bar{s}-s$ for $s\in\left[
\bar{s}-1,\bar{s}\right]  $, then by (\ref{369 894}),\vspace{-0.02in}%
\[
\int_{\bar{s}/2}^{\bar{s}-1}\left\vert \operatorname{Rc}\right\vert
^{2}\left(  \gamma\left(  s\right)  \right)  ds\leq\int_{0}^{\bar{s}}%
\zeta\left(  s\right)  ^{2}\left\vert \operatorname{Rc}\right\vert ^{2}\left(
\gamma\left(  s\right)  \right)  ds \leq \frac{n+\sqrt{1+4c}}{c}.\vspace{-0.02in}%
\]
Hence there is $z=\gamma(s_{0}), \, s_{0}\in\lbrack\bar{s}/2,\bar{s}-1]$, such
that $\left\vert \operatorname{Rc}\right\vert ^{2}(z)<\frac{4(n+\sqrt{1+4c}%
)}{\bar{s}c}$. Since $\left\vert S\right\vert \leq\sqrt{1+4c}$, $d(z,y)\leq
\int_{s_{0}}^{\bar{s}}\left\vert S\right\vert ds\leq\frac{\sqrt{1+4c}}%
{2}d(x,y)$ and $d(x,z)\leq(1+\frac{\sqrt{1+4c}}{2})d(x,y)$.

We have proved that given $x\in\mathcal{M}$, there exists
$\operatorname{const}<\infty$ such that for every $y\in\mathcal{M}\setminus
B(x,4)$ there exists $z\in\mathcal{M}$ with $d\left(  z,y\right)  \leq
\frac{\sqrt{1+4c}}{2}d\left(  x,y\right)  $ such that $\left\vert
\operatorname{Rc}\right\vert \left(  z\right)  \leq\frac{\operatorname{const}%
}{\sqrt{d\left(  x,z\right)  }}$. Hence $\liminf_{z\rightarrow\infty
}\left\vert \operatorname{Rc}\right\vert \left(  z\right)  =0$. This is a
result of Fern\'{a}ndez-Lopez and Garc\'{\i}a-R\'{\i}o \cite{Fernandex}. By
estimating the potential $f$, Wu \cite{PengWu} obtained $\liminf
_{z\rightarrow\infty}R\!\left(  z\right)  =0$.\smallskip

\textbf{Remark.} Define $\rho\left(  x,y,\bar{s}\right)  \!\doteqdot
\!\inf_{\gamma}\mathcal{J}\!\left(  \gamma\right)  $, where the inf is over
$\gamma:\left[  0,\bar{s}\right]  \rightarrow\mathcal{M}$ from $x$ to $y$. Let
$\gamma$ be a minimal $\phi$-geodesic. Since $\nabla_{S}S=\nabla\phi$, we have
$\tfrac{1}{2}\left.  \tfrac{d}{du}\right\vert _{u=0}\mathcal{J}\left(
\gamma_{u}\right)  =\left.  \left\langle U,S\right\rangle \right\vert
_{0}^{\bar{s}}$. Thus $\nabla\rho\left(  x,y,\bar{s}\right)  =2S\left(
\bar{s}\right)  =2\gamma^{\prime}\left(  \bar{s}\right)  $. By Lemma 3.1 of
\cite{LY}, we have $\frac{\partial\rho}{\partial\bar{s}}+\frac{1}{4}\left\vert
\nabla\rho\right\vert ^{2}=2\phi$ in the weak sense. Let $\left\{
e_{i}\right\}  _{i=1}^{n}$ be a parallel orthonormal frame along $\gamma_{0}$.
Summing $U=\frac{s}{\bar{s}}e_{i}$ over $i$ in the second variation formula,
we obtain\vspace{-0.02in}%
\[
\tfrac{1}{2}\Delta_{y}\rho\leq\frac{n}{\bar{s}}+\int_{0}^{\bar{s}}\left(
-\frac{s^{2}}{\bar{s}^{2}}\operatorname{Rc}\left(  S,S\right)  +\frac{s^{2}%
}{\bar{s}^{2}}\left\langle \nabla f,\nabla\phi\right\rangle +\frac{s^{2}}%
{\bar{s}^{2}}\Delta_{f}\phi\right)  ds.\vspace{-0.02in}%
\]
Since\vspace{-0.02in}%
\begin{align*}
-\!\int_{0}^{\bar{s}}\!\frac{s^{2}}{\bar{s}^{2}}\operatorname{Rc}\left(
S,S\right)  \!ds &  =\int_{0}^{\bar{s}}\frac{s^{2}}{\bar{s}^{2}}S\left(
S\left(  f\right)  \right)  ds-\int_{0}^{\bar{s}}\frac{s^{2}}{\bar{s}^{2}%
}\left\langle \nabla_{S}S,\nabla f\right\rangle ds\vspace{-0.02in}\\
&  =\left\langle S,\nabla f\right\rangle \!\left(  \bar{s}\right)  -\frac
{2}{\bar{s}}f\!\left(  y\right)  +\frac{2}{\bar{s}^{2}}\!\int_{0}^{\bar{s}%
}\!f\!\left(  \gamma\left(  s\right)  \right)  \!ds-\!\int_{0}^{\bar{s}%
}\!\frac{s^{2}}{\bar{s}^{2}}\!\left\langle \nabla\phi,\nabla f\right\rangle
\!ds,
\end{align*}
we have\vspace{-0.02in}%
\[
\tfrac{1}{2}\Delta_{y}\rho\left(  x,y,\bar{s}\right)  \leq\frac{n}{\bar{s}%
}+\int_{0}^{\bar{s}}\frac{s^{2}}{\bar{s}^{2}}\Delta_{f}\phi ds+\left\langle
S,\nabla f\right\rangle \left(  \bar{s}\right)  +\frac{2}{\bar{s}^{2}}\int
_{0}^{\bar{s}}\left(  f\left(  \gamma\left(  s\right)  \right)  -f\left(
y\right)  \right)  ds.\vspace{-0.02in}%
\]
Thus $\Phi=\bar{s}^{-\frac{n}{2}}e^{-\frac{\rho}{4}}$ satisfies in the weak
sense (note $\left\langle S,\nabla f\right\rangle \left(  \bar{s}\right)
=\frac{1}{2}\left\langle \nabla\rho,\nabla f\right\rangle $)\vspace{-0.02in}%
\[
\left(  \frac{\partial}{\partial\bar{s}}-\left(  \Delta_{y}\right)  _{f}%
+\frac{\phi}{2}\right)  \Phi\leq\frac{\Phi}{2}\left(  \int_{0}^{\bar{s}}%
\frac{s^{2}}{\bar{s}^{2}}\Delta_{f}\phi ds+\frac{2}{\bar{s}^{2}}\int_{0}%
^{\bar{s}}\left(  f\left(  \gamma\left(  s\right)  \right)  -f\left(
y\right)  \right)  ds\right)  .\vspace{-0.02in}%
\]

\textbf{Acknowledgment.} We wish to thank Jiaping Wang for helpful
discussions. We also wish to thank Bo Yang for correcting an error in the
previous version.

\end{document}